\documentclass[leqno,12pt]{amsart}
\usepackage{amssymb} 
\theoremstyle{plain} 
\newtheorem{theorem}{\indent\sc Theorem}[section] 

\newtheorem{proposition}[theorem]{\indent\sc Proposition}

\theoremstyle{definition} 
\newtheorem{definition}[theorem]{\indent\sc Definition}

\begin{document}

\title{Studies on the second member of the second Painlev\'e hierarchy \\}

\renewcommand{\thefootnote}{\fnsymbol{footnote}}
\footnote[0]{2000\textit{ Mathematics Subjet Classification}.
34M55; 34M45; 58F05; 32S65.}

\keywords{ 
Affine Weyl group, birational symmetry, coupled Painlev\'e system.}
\maketitle

\begin{abstract}
In this paper, we study the second member of the second Painlev\'e hierarchy $P_{II}^{(2)}$. We show that the birational transformations take this equation to the polynomial Hamiltonian system in dimension four, and this Hamiltonian system can be considered as a 1-parameter family of coupled Painlev\'e systems. This Hamiltonian is new. We also show that this system admits extended affine Weyl group symmetry of type $A_1^{(1)}$, and can be recovered by its holomorphy conditions. We also study a fifth-order ordinary differential equation satisfied by this Hamiltonian. After we transform this equation into a system of the first-order ordinary differential equations of polynomial type in dimension five by birational transformations, we give its symmetry and holomorphy conditions.
\end{abstract}

\section{Introduction}

In this paper, we study the second member of the second Painlev\'e hierarchy \cite{Clarkson,Mazzocco,Joshi} explicitly given by
\begin{equation}\label{eq:1}
P_{II}^{(2)}:\frac{d^4u}{dt^4}=10u \left(\frac{du}{dt} \right)^2+10u^2 \frac{d^2u}{dt^2}-6u^5+tu+\alpha_2 \quad (\alpha_2 \in {\Bbb C}).
\end{equation}
It was found by Martynov in 1973 (see \cite{Martynov}) and rediscovered independently by Ablowitz and Segur in 1977 (see \cite{Ablowitz}) and Flaschka and Newell in 1980 (see \cite{Flaschka}). It is known that this system can be obtained by self-similar reduction of the Modified KdV5 equation
\begin{equation}\label{Kdv}
g_{wwwww}=10(g^2-2R(s))g_{www}+40gg_{w}g_{ww}+10{g_{w}}^3-30(g^2-2R(s))^2 g_{w}+g_{s},
\end{equation}
where $R(s)$ is an arbitrary locally analytic function of $s$.

We note that this equation appears as the equation F-XVII in Cosgrove's classification of the fourth-order ordinary differential equations in the polynomial class having the Painlev\'e property (see \cite{Cosgrove}).

At first, we show that the birational transformations (see Section 2) take the equation \eqref{eq:1} to a polynomial Hamiltonian system in dimension four. We make this polynomial Hamiltonian from the viewpoint of accessible singularity and local index (see Section 5).

This Hamiltonian system can be considered as a 1-parameter family of coupled Painlev\'e systems in dimension four. This Hamiltonian is new.

It is known that for the equation \eqref{Kdv} M. Mazzocco and M. Y. Mo found a canonical variables $(q_i,p_i)_{i=1}^2$ and obtained a polynomial Hamiltonian through the monodromy preserving deformation equation of the second-order ordinary differential equation (see \cite{Mazzocco}). The degree of this Hamiltonian with respect to $q_1,p_1,q_2,p_2$ is 4. On the other hand, one of our Hamiltonian is 3.

We also study its symmetry and holomorphy conditions. We show that this system admits extended affine Weyl group symmetry of type $A_1^{(1)}$ as the group of its B{\"a}cklund transformations. These B{\"a}cklund transformations satisfy
\begin{equation}
s_i(g)=g+\frac{\alpha_i}{f_i}\{f_i,g\}+\frac{1}{2!} \left(\frac{\alpha_i}{f_i} \right)^2 \{f_i,\{f_i,g\} \}+\cdots \quad (g \in {\Bbb C}(t)[q_1,p_1,q_2,p_2]),
\end{equation}
where poisson bracket $\{,\}$ satisfies the relations:
$$
\{p_1,q_1\}=\{p_2,q_2\}=1, \quad the \ others \ are \ 0.
$$
Since these B{\"a}cklund transformations have Lie theoretic origin, similarity reduction of a Drinfeld-Sokolov hierarchy admits such a B{\"a}cklund symmetry.

These properties of its symmetry and holomorphy conditions are new.

We find a rational solution of this Hamiltonian system given by
\begin{equation}
(q_1,p_1,q_2,p_2;\alpha_2)=\left(0,-\frac{t}{2},0,0;0\right)
\end{equation}
as a seed solution. Applying the translation operators, we can obtain an infinite series of the rational solutions.

Finally, we also study a fifth-order ordinary differential equation satisfied by this Hamiltonian. We show that the birational transformations take this equation to a system of the first-order ordinary differential equations of polynomial type in dimension five. For this system, we give its symmetry and holomorphy conditions.

It is still the following open questions on the system \eqref{eq:3}:
\begin{enumerate}
\item Relation between Mazzocco's Hamiltonian (see \cite{Mazzocco}) with \eqref{eq:4}.
\item Special polynomials.
\item Lax pair.
\item Reduction involving the $A_1^{(1)}$-symmetry.
\item Bilinear form.
\item Classification of the rational solutions or special solutions.
\item Irreducibility.
\item Discrete version.
\item q-Analogue.
\item Quantum version.
\item The remained members of $P_{II}^{(n)}$.
\end{enumerate}

\section{Polynomial Hamiltonian}
\begin{theorem}\label{th1.1}
The birational transformations
\begin{equation}\label{eq:2}
  \left\{
  \begin{aligned}
   q_1 =&u,\\
   p_1 =&\frac{d^3u}{dt^3}+\left(\frac{du}{dt} \right)^2-\frac{t}{2}+\left(3u^3-6u\frac{du}{dt}-2\frac{d^2u}{dt^2} \right)u,\\
   q_2 =&\frac{d^2u}{dt^2}-2u\frac{du}{dt},\\
   p_2 =&\frac{du}{dt}-u^2
   \end{aligned}
  \right. 
\end{equation}
take the system to the Hamiltonian system
\begin{equation}\label{eq:3}
  \left\{
  \begin{aligned}
   \frac{dq_1}{dt} =&\frac{\partial H}{\partial p_1}=q_1^2+p_2,\\
   \frac{dp_1}{dt} =&-\frac{\partial H}{\partial q_1}=-2q_1p_1+\alpha_2-\frac{1}{2},\\
   \frac{dq_2}{dt} =&\frac{\partial H}{\partial p_2}=-3p_2^2+p_1+\frac{t}{2},\\
   \frac{dp_2}{dt} =&-\frac{\partial H}{\partial q_2}=q_2
   \end{aligned}
  \right. 
\end{equation}
with the polynomial Hamiltonian \rm{(cf. \cite{Mazzocco}) \rm}
\begin{align}\label{eq:4}
\begin{split}
H=&K(q_1,p_1;\alpha_2)+H_I(q_2,p_2,t)+p_1p_2\\
=&q_1^2p_1+\left(\frac{1}{2}-\alpha_2 \right)q_1-p_2^3+\frac{t}{2}p_2-\frac{q_2^2}{2}+p_1p_2.
\end{split}
\end{align}
\end{theorem}
The symbols $K(x,y;\alpha)$ and $H_I(z,w,t)$ denote
\begin{align}
\begin{split}
K(x,y;\alpha)=&x^2y+\left(\frac{1}{2}-\alpha \right)x,\\
H_I(z,w,t)=&-w^3+\frac{t}{2}w-\frac{z^2}{2}.
\end{split}
\end{align}
The system with the Hamiltonian $K(x,y;\alpha)$ has itself as its first integral, and  $H_I(z,w,t)$ denotes the Painlev\'e I Hamiltonian.

This system is a 1-parameter family of coupled Painlev\'e systems. This Hamiltonian is new  (cf. \cite{Mazzocco}).

Before we will prove Theorem \ref{th1.1}, we review the notion of accessible singularity and local index.

\section{Accessible singularity and local index}
Let us review the notion of {\it accessible singularity}. Let $B$ be a connected open domain in $\Bbb C$ and $\pi : {\mathcal W} \longrightarrow B$ a smooth proper holomorphic map. We assume that ${\mathcal H} \subset {\mathcal W}$ is a normal crossing divisor which is flat over $B$. Let us consider a rational vector field $\tilde v$ on $\mathcal W$ satisfying the condition
\begin{equation*}
\tilde v \in H^0({\mathcal W},\Theta_{\mathcal W}(-\log{\mathcal H})({\mathcal H})).
\end{equation*}
Fixing $t_0 \in B$ and $P \in {\mathcal W}_{t_0}$, we can take a local coordinate system $(x_1,\ldots ,x_n)$ of ${\mathcal W}_{t_0}$ centered at $P$ such that ${\mathcal H}_{\rm smooth \rm}$ can be defined by the local equation $x_1=0$.
Since $\tilde v \in H^0({\mathcal W},\Theta_{\mathcal W}(-\log{\mathcal H})({\mathcal H}))$, we can write down the vector field $\tilde v$ near $P=(0,\ldots ,0,t_0)$ as follows:
\begin{equation*}
\tilde v= \frac{\partial}{\partial t}+g_1 
\frac{\partial}{\partial x_1}+\frac{g_2}{x_1} 
\frac{\partial}{\partial x_2}+\cdots+\frac{g_n}{x_1} 
\frac{\partial}{\partial x_n}.
\end{equation*}
This vector field defines the following system of differential equations
\begin{equation}\label{39}
\frac{dx_1}{dt}=g_1(x_1,\ldots,x_n,t),\ \frac{dx_2}{dt}=\frac{g_2(x_1,\ldots,x_n,t)}{x_1},\cdots, \frac{dx_n}{dt}=\frac{g_n(x_1,\ldots,x_n,t)}{x_1}.
\end{equation}
Here $g_i(x_1,\ldots,x_n,t), \ i=1,2,\dots ,n,$ are holomorphic functions defined near $P=(0,\dots ,0,t_0).$

\begin{definition}\label{Def1}
With the above notation, assume that the rational vector field $\tilde v$ on $\mathcal W$ satisfies the condition
$$
(A) \quad \tilde v \in H^0({\mathcal W},\Theta_{\mathcal W}(-\log{\mathcal H})({\mathcal H})).
$$
We say that $\tilde v$ has an {\it accessible singularity} at $P=(0,\dots ,0,t_0)$ if
$$
x_1=0 \ {\rm and \rm} \ g_i(0,\ldots,0,t_0)=0 \ {\rm for \rm} \ {\rm every \rm} \ i, \ 2 \leq i \leq n.
$$
\end{definition}

If $P \in {\mathcal H}_{{\rm smooth \rm}}$ is not an accessible singularity, all solutions of the ordinary differential equation passing through $P$ are vertical solutions, that is, the solutions are contained in the fiber ${\mathcal W}_{t_0}$ over $t=t_0$. If $P \in {\mathcal H}_{\rm smooth \rm}$ is an accessible singularity, there may be a solution of \eqref{39} which passes through $P$ and goes into the interior ${\mathcal W}-{\mathcal H}$ of ${\mathcal W}$.

Here we review the notion of {\it local index}. Let $v$ be an algebraic vector field with an accessible singular point $\overrightarrow{p}=(0,\ldots,0)$ and $(x_1,\ldots,x_n)$ be a coordinate system in a neighborhood centered at $\overrightarrow{p}$. Assume that the system associated with $v$ near $\overrightarrow{p}$ can be written as
\begin{align}\label{b}
\begin{split}
\frac{d}{dt}\begin{pmatrix}
             x_1 \\
             x_2 \\
             \vdots\\
             x_{n-1} \\
             x_n
             \end{pmatrix}=\frac{1}{x_1}\left\{\begin{bmatrix}
             a_{11} & 0 & 0 & \hdots & 0 \\
             a_{21} & a_{22} & 0 &  \hdots & 0 \\
             \vdots & \vdots & \ddots & 0 & 0 \\
             a_{(n-1)1} & a_{(n-1)2} & \hdots & a_{(n-1)(n-1)} & 0 \\
             a_{n1} & a_{n2} & \hdots & a_{n(n-1)} & a_{nn}
             \end{bmatrix}\begin{pmatrix}
             x_1 \\
             x_2 \\
             \vdots\\
             x_{n-1} \\
             x_n
             \end{pmatrix}+\begin{pmatrix}
             x_1h_1(x_1,\ldots,x_n,t) \\
             h_2(x_1,\ldots,x_n,t) \\
             \vdots\\
             h_{n-1}(x_1,\ldots,x_n,t) \\
             h_n(x_1,\ldots,x_n,t)
             \end{pmatrix}\right\},\\
              (h_i \in {\Bbb C}(t)[x_1,\ldots,x_n], \ a_{ij} \in {\Bbb C}(t))
             \end{split}
             \end{align}
where $h_1$ is a polynomial which vanishes at $\overrightarrow{p}$ and $h_i$, $i=2,3,\ldots,n$ are polynomials of order at least 2 in $x_1,x_2,\ldots,x_n$, We call ordered set of the eigenvalues $(a_{11},a_{22},\cdots,a_{nn})$ {\it local index} at $\overrightarrow{p}$.

We are interested in the case with local index
\begin{equation}\label{integer}
(1,a_{22}/a_{11},\ldots,a_{nn}/a_{11}) \in {\Bbb Z}^{n}.
\end{equation}
These properties suggest the possibilities that $a_1$ is the residue of the formal Laurent series:
\begin{equation}
y_1(t)=\frac{a_{11}}{(t-t_0)}+b_1+b_2(t-t_0)+\cdots+b_n(t-t_0)^{n-1}+\cdots \quad (b_i \in {\Bbb C}),
\end{equation}
and the ratio $(1,a_{22}/a_{11},\ldots,a_{nn}/a_{11})$ is resonance data of the formal Laurent series of each $y_i(t) \ (i=2,\ldots,n)$, where $(y_1,\ldots,y_n)$ is original coordinate system satisfying $(x_1,\ldots,x_n)=(f_1(y_1,\ldots,y_n),\ldots,f_n(y_1,\ldots,y_n)), \ f_i(y_1,\ldots,y_n) \in {\Bbb C}(t)(y_1,\ldots,y_n)$.

If each component of $(1,a_{22}/a_{11},\ldots,a_{nn}/a_{11})$ has the same sign, we may resolve the accessible singularity by blowing-up finitely many times. However, when different signs appear, we may need to both blow up and blow down.

The $\alpha$-test,
\begin{equation}\label{poiuy}
t=t_0+\alpha T, \quad x_i=\alpha X_i, \quad \alpha \rightarrow 0,
\end{equation}
yields the following reduced system:
\begin{align}\label{ppppppp}
\begin{split}
\frac{d}{dT}\begin{pmatrix}
             X_1 \\
             X_2 \\
             \vdots\\
             X_{n-1} \\
             X_n
             \end{pmatrix}=\frac{1}{X_1}\begin{bmatrix}
             a_{11}(t_0) & 0 & 0 & \hdots & 0 \\
             a_{21}(t_0) & a_{22}(t_0) & 0 &  \hdots & 0 \\
             \vdots & \vdots & \ddots & 0 & 0 \\
             a_{(n-1)1}(t_0) & a_{(n-1)2}(t_0) & \hdots & a_{(n-1)(n-1)}(t_0) & 0 \\
             a_{n1}(t_0) & a_{n2}(t_0) & \hdots & a_{n(n-1)}(t_0) & a_{nn}(t_0)
             \end{bmatrix}\begin{pmatrix}
             X_1 \\
             X_2 \\
             \vdots\\
             X_{n-1} \\
             X_n
             \end{pmatrix},
             \end{split}
             \end{align}
where $a_{ij}(t_0) \in {\Bbb C}$. Fixing $t=t_0$, this system is the system of the first order ordinary differential equation with constant coefficient. Let us solve this system. At first, we solve the first equation:
\begin{equation}
X_1(T)=a_{11}(t_0)T+C_1 \quad (C_1 \in {\Bbb C}).
\end{equation}
Substituting this into the second equation in \eqref{ppppppp}, we can obtain the first order linear ordinary differential equation:
\begin{equation}
\frac{dX_2}{dT}=\frac{a_{22}(t_0) X_2}{a_{11}(t_0)T+C_1}+a_{21}(t_0).
\end{equation}
By variation of constant, in the case of $a_{11}(t_0) \not= a_{22}(t_0)$ we can solve explicitly:
\begin{equation}
X_2(T)=C_2(a_{11}(t_0)T+C_1)^{\frac{a_{22}(t_0)}{a_{11}(t_0)}}+\frac{a_{21}(t_0)(a_{11}(t_0)T+C_1)}{a_{11}(t_0)-a_{22}(t_0)} \quad (C_2 \in {\Bbb C}).
\end{equation}
This solution is a single-valued solution if and only if
$$
\frac{a_{22}(t_0)}{a_{11}(t_0)} \in {\Bbb Z}.
$$
In the case of $a_{11}(t_0)=a_{22}(t_0)$ we can solve explicitly:
\begin{equation}
X_2(T)=C_2(a_{11}(t_0)T+C_1)+\frac{a_{21}(t_0)(a_{11}(t_0)T+C_1){\rm Log}(a_{11}(t_0)T+C_1)}{a_{11}(t_0)} \quad (C_2 \in {\Bbb C}).
\end{equation}
This solution is a single-valued solution if and only if
$$
a_{21}(t_0)=0.
$$
Of course, $\frac{a_{22}(t_0)}{a_{11}(t_0)}=1 \in {\Bbb Z}$.
In the same way, we can obtain the solutions for each variables $(X_3,\ldots,X_n)$. The conditions $\frac{a_{jj}(t)}{a_{11}(t)} \in {\Bbb Z}, \ (j=2,3,\ldots,n)$ are necessary condition in order to have the Painlev\'e property.

\section{The case of the second Painlev\'e system}

In this section, we review the case of the second Painlev\'e system:
\begin{equation}\label{PII}
\frac{d^2u}{dt^2}=2u^3+tu+\alpha \quad (\alpha \in {\Bbb C}).
\end{equation}

Let us make its polynomial Hamiltonian from the viewpoint of accessible singularity and local index.

{\bf Step 0:} We make a change of variables.
\begin{equation}
x=u, \quad y=\frac{du}{dt}.
\end{equation}

{\bf Step 1:} We make a change of variables.
\begin{equation}
x_1=\frac{1}{x}, \quad y_1=\frac{y}{x^2}.
\end{equation}
In this coordinate system, we see that this system has two accessible singular points:
\begin{equation}
(x_1,y_1)=\left\{(0,1),(0,-1) \right\}.
\end{equation}

Around the point $(x_1,y_1)=(0,1)$, we can rewrite the system as follows.

{\bf Step 2:} We make a change of variables.
\begin{equation}
x_2=x_1, \quad y_2=y_1-1.
\end{equation}
In this coordinate system, we can rewrite the system satisfying the condition \eqref{b}:
\begin{align*}
\frac{d}{dt}\begin{pmatrix}
             x_2 \\
             y_2 
             \end{pmatrix}&=\frac{1}{x_2}\left\{\begin{pmatrix}
             -1 & 0   \\
             0 & -4 
             \end{pmatrix}\begin{pmatrix}
             x_2 \\
             y_2 
             \end{pmatrix}+\cdots\right\},
             \end{align*}
and we can obtain the local index $(-1,-4)$ at the point $\{(x_2,y_2)=(0,0)\}$. The ratio of the local index at the point $\{(x_2,y_2)=(0,0)\}$ is a positive integer.

We aim to obtain the local index $(-1,-2)$ by successive blowing-up procedures.

{\bf Step 3:} We blow up at the point $\{(x_2,y_2)=(0,0)\}$.
\begin{equation}
x_3=x_2, \quad y_3=\frac{y_2}{x_2}.
\end{equation}

{\bf Step 4:} We blow up at the point $\{(x_3,y_3)=(0,0)\}$.
\begin{equation}
x_4=x_3, \quad y_4=\frac{y_3}{x_3}.
\end{equation}
In this coordinate system, we see that this system has the following accessible singular point:
\begin{equation}
(x_4,y_4)=(0,t/2).
\end{equation}

{\bf Step 5:} We make a change of variables.
\begin{equation}
x_5=x_4, \quad y_5=y_4-t/2.
\end{equation}
In this coordinate system, we can rewrite the system as follows:
\begin{align*}
\frac{d}{dt}\begin{pmatrix}
             x_5 \\
             y_5 
             \end{pmatrix}&=\frac{1}{x_5}\left\{\begin{pmatrix}
             -1 & 0  \\
             \alpha-1/2 & -2
             \end{pmatrix}\begin{pmatrix}
             x_5 \\
             y_5 
             \end{pmatrix}+\cdots\right\},
             \end{align*}
and we can obtain the local index $(-1,-2)$. Here, the relation between $(x_5,y_5)$ and $(x,y)$ is given by
\begin{equation*}
  \left\{
  \begin{aligned}
   x_5 &=\frac{1}{x},\\
   y_5 &=y-x^2-\frac{t}{2}.
   \end{aligned}
  \right. 
\end{equation*}

Finally, we can choose canonical variables $(q,p)$.

{\bf Step 9:} We make a change of variables.
\begin{equation}
q=\frac{1}{x_5}, \quad p=y_5,
\end{equation}
and we can obtain the system
\begin{equation*}
  \left\{
  \begin{aligned}
   \frac{dq}{dt} &=q^2+p+\frac{t}{2},\\
   \frac{dp}{dt} &=-2qp+\alpha-\frac{1}{2}
   \end{aligned}
  \right. 
\end{equation*}
with the polynomial Hamiltonian $H_{II}$:
\begin{equation}
H_{II}=q^2 p+\frac{1}{2}p^2+\frac{t}{2}p-\left(\alpha-\frac{1}{2} \right)q.
\end{equation}
We remark that we can discuss the case of the accessible singular point $(x_1,y_1)=(0,-1)$ in the same way as in the case of $(x_1,y_1)=(0,1)$.

\section{Proof of theorem \ref{th1.1}}

By the same way of the second Painlev\'e system, we can prove Theorem \ref{th1.1}.

At first, we rewrite the equation \eqref{eq:2} to the system of the first order ordinary differential equations.

{\bf Step 0:} We make a change of variables.
\begin{equation}
x=u, \quad y=\frac{du}{dt}, \quad z=\frac{d^2u}{dt^2}, \quad w=\frac{d^3u}{dt^3}.
\end{equation}

{\bf Step 1:} We make a change of variables.
\begin{equation}
x_1=\frac{1}{x}, \quad y_1=\frac{y}{x^2}, \quad z_1=\frac{z}{x^3}, \quad w_1=\frac{w}{x^4}.
\end{equation}
In this coordinate system, we see that this system has four accessible singular points:
\begin{equation}
(x_1,y_1,z_1,w_1)=\left\{(0,1,2,6),(0,-1,2,-6),\left(0,\frac{1}{2},\frac{1}{2},\frac{3}{4} \right),\left(0,-\frac{1}{2},\frac{1}{2},-\frac{3}{4} \right) \right\}.
\end{equation}

Around the point $(x_1,y_1,z_1,w_1)=(0,-1,2,-6)$, we can rewrite the system as follows.
{\bf Step 2:} We make a change of variables.
\begin{equation}
x_2=x_1, \quad y_2=y_1+1, \quad z_2=z_1-2, \quad w_2=w_1+6.
\end{equation}
In this coordinate system, we can rewrite the system satisfying the condition \eqref{b}:
\begin{align*}
\frac{d}{dt}\begin{pmatrix}
             x_2 \\
             y_2 \\
             z_2 \\
             w_2
             \end{pmatrix}&=\frac{1}{x_2}\left\{\begin{pmatrix}
             1 & 0 & 0 & 0  \\
             0 & 4 & 1 & 0 \\
             0 & -6 & 3 & 1 \\
             0 & 4 & 10 & 4
             \end{pmatrix}\begin{pmatrix}
             x_2 \\
             y_2 \\
             z_2 \\
             w_2 
             \end{pmatrix}+\cdots\right\}.
             \end{align*}
To the above system, we make the linear transformation:
\begin{align*}
\begin{pmatrix}
             X_2 \\
             Y_2 \\
             Z_2 \\
             W_2
             \end{pmatrix}&=\begin{pmatrix}
             1 & 0 & 0 & 0  \\
             0 & 1 & 1 & 1 \\
             0 & -2 & -1 & 2 \\
             0 & 8 & 6 & 12
             \end{pmatrix}\begin{pmatrix}
             x_2 \\
             y_2 \\
             z_2 \\
             w_2 
             \end{pmatrix}
             \end{align*}
to arrive at
\begin{align*}
\frac{d}{dt}\begin{pmatrix}
             X_2 \\
             Y_2 \\
             Z_2 \\
             W_2
             \end{pmatrix}&=\frac{1}{X_2}\left\{\begin{pmatrix}
             1 & 0 & 0 & 0  \\
             0 & 2 & 0 & 0 \\
             0 & 0 & 3 & 0 \\
             0 & 0 & 0 & 6
             \end{pmatrix}\begin{pmatrix}
             X_2 \\
             Y_2 \\
             Z_2 \\
             W_2 
             \end{pmatrix}+\cdots\right\},
             \end{align*}
and we can obtain the local index $(1,2,3,6)$ at the point $\{(X_2,Y_2,Z_2,W_2)=(0,0,0,0)\}$. The continued ratio of the local index at the point $\{(X_2,Y_2,Z_2,W_2)=(0,0,0,0)\}$ are all positive integers
\begin{equation}
\left(\frac{2}{1},\frac{3}{1},\frac{6}{1} \right)=(2,3,6).
\end{equation}
This is the reason why we choose this accessible singular point.

We aim to obtain the local index $(1,0,0,2)$ by successive blowing-up procedures.

{\bf Step 3:} We blow up at the point $\{(x_2,y_2,z_2,w_2)=(0,0,0,0)\}$.
\begin{equation}
x_3=x_2, \quad y_3=\frac{y_2}{x_2}, \quad z_3=\frac{z_2}{x_2}, \quad w_3=\frac{w_2}{x_2}.
\end{equation}

{\bf Step 4:} We blow up at the point $\{(x_3,y_3,z_3,w_3)=(0,0,0,0)\}$.
\begin{equation}
x_4=x_3, \quad y_4=\frac{y_3}{x_3}, \quad z_4=\frac{z_3}{x_3}, \quad w_4=\frac{w_3}{x_3}.
\end{equation}
In this coordinate system, we see that this system has the following accessible singular locus:
\begin{equation}
(x_4,y_4,z_4,w_4)=(0,y_4,-2y_4,8y_4).
\end{equation}

{\bf Step 5:} We blow up along the curve $\{(x_4,y_4,z_4,w_4)=(0,y_4,-2y_4,8y_4)\}$.
\begin{equation}
x_5=x_4, \quad y_5=y_4, \quad z_5=\frac{z_4+2y_4}{x_4}, \quad w_5=\frac{w_4-8y_4}{x_4}.
\end{equation}
In this coordinate system, we see that this system has the following accessible singular locus:
\begin{equation}
(x_5,y_5,z_5,w_5)=(0,y_5,z_5,-2z_5).
\end{equation}

{\bf Step 6:} We blow up along the surface $\{(x_5,y_5,z_5,w_5)=(0,y_5,z_5,-2z_5)\}$.
\begin{equation}
x_6=x_5, \quad y_6=y_5, \quad z_6=z_5, \quad w_6=\frac{w_5+2z_5}{x_5}.
\end{equation}
In this coordinate system, we see that this system has the following accessible singular locus:
\begin{equation}
(x_6,y_6,z_6,w_6)=\left(0,y_6,z_6,y_6^2-\frac{t}{2} \right).
\end{equation}

{\bf Step 7:} We make a change of variables.
\begin{equation}
x_7=x_6, \quad y_7=y_6, \quad z_7=z_6, \quad w_7=w_6-y_6^2+\frac{t}{2}.
\end{equation}
In this coordinate system, we can rewrite the system as follows:
\begin{align*}
\frac{d}{dt}\begin{pmatrix}
             x_7 \\
             y_7 \\
             z_7 \\
             w_7
             \end{pmatrix}&=\frac{1}{x_7}\left\{\begin{pmatrix}
             1 & 0 & 0 & 0  \\
             0 & 0 & 0 & 0 \\
             -\frac{t}{2} & 0 & 0 & 0 \\
             \alpha_2+\frac{1}{2} & 0 & 0 & 2
             \end{pmatrix}\begin{pmatrix}
             x_7 \\
             y_7 \\
             z_7 \\
             w_7 
             \end{pmatrix}+\cdots\right\},
             \end{align*}
and we can obtain the local index $(1,0,0,2)$. Here, the relation between $(x_7,y_7,z_7,w_7)$ and $(x,y,z,w)$ is given by
\begin{equation*}
  \left\{
  \begin{aligned}
   x_7 &=\frac{1}{x},\\
   y_7 &=x^2+y,\\
   z_7 &=z+2xy,\\
   w_7 &=w+\frac{t}{2}-3x^4-6x^2 y-y^2+2xz.
   \end{aligned}
  \right. 
\end{equation*}

{\bf Step 8:} We make a change of variables.
\begin{equation}
x_8=\frac{1}{x_7}, \quad y_8=y_7, \quad z_8=z_7, \quad w_8=w_7.
\end{equation}
In this coordinate system, we can rewrite the system as follows:
\begin{equation*}
  \left\{
  \begin{aligned}
   \frac{dx_8}{dt} &=-x_8^2+y_8,\\
    \frac{dy_8}{dt} &=z_8,\\
    \frac{dz_8}{dt} &=3y_8^2+w_8-\frac{t}{2},\\
    \frac{dw_8}{dt} &=2x_8 w_8+\alpha_2+\frac{1}{2}.
   \end{aligned}
  \right. 
\end{equation*}
Finally, we can choose canonical variables $(q_1,p_1,q_2,p_2)$.

{\bf Step 9:} We make a change of variables.
\begin{equation}
q_1=-x_8, \quad p_1=-w_8, \quad q_2=-z_8, \quad p_2=-y_8,
\end{equation}
and we can obtain the system \eqref{eq:3} with the polynomial Hamiltonian \eqref{eq:4}.

Thus, we have completed the proof of Theorem \ref{th1.1} \qed.

We note on the remaining accessible singular points.

Around the point $(x_1,y_1,z_1,w_1)=(0,1,2,6)$, we can rewrite the system as follows.

{\bf Step 2:} We make a change of variables.
\begin{equation}
x_2=x_1, \quad y_2=y_1-1, \quad z_2=z_1-2, \quad w_2=w_1-6.
\end{equation}
In this coordinate system, we can rewrite the system satisfying the condition \eqref{b}:
\begin{align*}
\frac{d}{dt}\begin{pmatrix}
             x_2 \\
             y_2 \\
             z_2 \\
             w_2
             \end{pmatrix}&=\frac{1}{x_2}\left\{\begin{pmatrix}
             -1 & 0 & 0 & 0  \\
             0 & -4 & 1 & 0 \\
             0 & -6 & -3 & 1 \\
             0 & -4 & 10 & -4
             \end{pmatrix}\begin{pmatrix}
             x_2 \\
             y_2 \\
             z_2 \\
             w_2 
             \end{pmatrix}+\cdots\right\}.
             \end{align*}
and we can obtain the local index $(-1,-2,-3,-6)$ at the point $\{(x_2,y_2,z_2,w_2)=(0,0,0,0)\}$. The continued ratio of the local index at the point $\{(x_2,y_2,z_2,w_2)=(0,0,0,0)\}$ are all positive integers
\begin{equation}
\left(\frac{-2}{-1},\frac{-3}{-1},\frac{-6}{-1} \right)=(2,3,6).
\end{equation}
We remark that we can discuss this case in the same way as in the case of $(x_1,y_1,z_1,w_1)=(0,-1,2,-6)$.

Around the point $(x_1,y_1,z_1,w_1)=(0,1/2,1/2,3/4)$, we can rewrite the system as follows.

{\bf Step 2:} We make a change of variables.
\begin{equation}
x_2=x_1, \quad y_2=y_1-1/2, \quad z_2=z_1-1/2, \quad w_2=w_1-3/4.
\end{equation}
In this coordinate system, we can rewrite the system satisfying the condition \eqref{b}:
\begin{align*}
\frac{d}{dt}\begin{pmatrix}
             x_2 \\
             y_2 \\
             z_2 \\
             w_2
             \end{pmatrix}&=\frac{1}{x_2}\left\{\begin{pmatrix}
             -1/2 & 0 & 0 & 0  \\
             0 & -2 & 1 & 0 \\
             0 & -3/2 & -3/2 & 1 \\
             0 & 7 & 10 & -2
             \end{pmatrix}\begin{pmatrix}
             x_2 \\
             y_2 \\
             z_2 \\
             w_2 
             \end{pmatrix}+\cdots\right\}.
             \end{align*}
and we can obtain the local index $(-1/2,3/2,-3,-4)$ at the point $\{(x_2,y_2,z_2,w_2)=(0,0,0,0)\}$. The continued ratio of the local index at the point $\{(x_2,y_2,z_2,w_2)=(0,0,0,0)\}$ are
\begin{equation}
\left(\frac{3/2}{-1/2},\frac{-3}{-1/2},\frac{-4}{-1/2} \right)=(-3,6,8).
\end{equation}
In this case, the local index involves a negative integer. So, we need to blow down.

Around the point $(x_1,y_1,z_1,w_1)=(0,-1/2,1/2,-3/4)$, we can rewrite the system as follows.

{\bf Step 2:} We make a change of variables.
\begin{equation}
x_2=x_1, \quad y_2=y_1+1/2, \quad z_2=z_1-1/2, \quad w_2=w_1+3/4.
\end{equation}
In this coordinate system, we can rewrite the system satisfying the condition \eqref{b}:
\begin{align*}
\frac{d}{dt}\begin{pmatrix}
             x_2 \\
             y_2 \\
             z_2 \\
             w_2
             \end{pmatrix}&=\frac{1}{x_2}\left\{\begin{pmatrix}
             1/2 & 0 & 0 & 0  \\
             0 & 2 & 1 & 0 \\
             0 & -3/2 & 3/2 & 1 \\
             0 & -7 & 10 & 2
             \end{pmatrix}\begin{pmatrix}
             x_2 \\
             y_2 \\
             z_2 \\
             w_2 
             \end{pmatrix}+\cdots\right\}.
             \end{align*}
and we can obtain the local index $(1/2,-3/2,3,4)$ at the point $\{(x_2,y_2,z_2,w_2)=(0,0,0,0)\}$. The continued ratio of the local index at the point $\{(x_2,y_2,z_2,w_2)=(0,0,0,0)\}$ are
\begin{equation}
\left(\frac{-3/2}{1/2},\frac{3}{1/2},\frac{4}{1/2} \right)=(-3,6,8).
\end{equation}
In this case, the local index involves a negative integer. So, we need to blow down.

\section{Symmetry and holomorphy conditions}
In this section, we study the symmetry and holomorphy conditions of the system \eqref{eq:3}. These symmetries, holomorphy conditions and invariant divisors are new.
\begin{theorem}\label{pro:3}
Let us consider a polynomial Hamiltonian system with Hamiltonian $H \in {\Bbb C}(t)[q_1,p_1,q_2,p_2]$. We assume that

$(A1)$ $deg(H)=5$ with respect to $q_1,p_1,q_2,p_2$.

$(A2)$ This system becomes again a polynomial Hamiltonian system in each coordinate system $r_i \ (i=0,1)${\rm : \rm}
\begin{align}\label{holo1}
\begin{split}
r_1:(x_1,y_1,z_1,w_1)=&\left(\frac{1}{q_1},-\left(q_1p_1+\frac{1}{2}-\alpha_2 \right)q_1,q_2,p_2 \right),\\
r_2:(x_2,y_2,z_2,w_2)=&(\frac{1}{q_1},-\left((p_1-2p_2^2+t+4q_1(q_1 p_2+q_2))q_1+1/2+\alpha_2 \right)q_1,\\
&q_2+4q_1(q_1^2+p_2),p_2+2q_1^2).
\end{split}
\end{align}
Then such a system coincides with the system \eqref{eq:3} with the polynomial Hamiltonian \eqref{eq:4}.
\end{theorem}
We note that the condition $(A2)$ should be read that
\begin{align*}
&r_0(H), \quad r_1 \left(H-q_1 \right)
\end{align*}
are polynomials with respect to $x_i,y_i,z_i,w_i$.

\begin{theorem}\label{th:2}
The system \eqref{eq:3} admits extended affine Weyl group symmetry of type $A_1^{(1)}$ as the group of its B{\"a}cklund transformations, whose generators $s_0,s_1,{\pi}$ defined as follows$:$ with {\it the notation} $(*):=(q_1,p_1,q_2,p_2,t;\alpha_2)$\rm{: \rm}
\begin{align*}
s_0:(*) \rightarrow &\left(q_1+\frac{\frac{1}{2}-\alpha_2}{p_1},p_1,q_2,p_2,t;1-\alpha_2 \right),\\
s_1:(*) \rightarrow &(q_1+\frac{2\alpha_2+1}{2(p_1+t-2p_2^2+4q_1(q_2+q_1p_2))},\\
&p_1-\frac{2(2\alpha_2+1)(q_2+2q_1p_2)}{p_1+t-2p_2^2+4q_1(q_2+q_1p_2)}+\frac{(2\alpha_2+1)^2(p_2+2q_1^2)}{(p_1+t-2p_2^2+4q_1(q_2+q_1p_2))^2},\\
&q_2-\frac{2(2\alpha_2+1)(p_2-q_1^2)}{p_1+t-2p_2^2+4q_1(q_2+q_1p_2)}+\frac{3(2\alpha_2+1)^2 q_1}{(p_1+t-2p_2^2+4q_1(q_2+q_1p_2))^2}\\
&+\frac{(2\alpha_2+1)^3}{2(p_1+t-2p_2^2+4q_1(q_2+q_1p_2))^3},\\
&p_2-\frac{2(2\alpha_2+1)q_1}{p_1+t-2p_2^2+4q_1(q_2+q_1p_2)}-\frac{(2\alpha_2+1)^2}{(p_1+t-2p_2^2+4q_1(q_2+q_1p_2))^2},\\
&t;-1-\alpha_2),\\
\pi:(*) &\rightarrow (-q_1,-(p_1+t-2p_2^2+4q_1(q_2+q_1p_2)),-(q_2+4q_1(q_1^2+p_2)),\\
&-(p_2+2q_1^2),t;-\alpha_2).
\end{align*}
\end{theorem}

\begin{proposition}\label{pro:2}
Let us define the following translation operators
\begin{equation}
T_1:=\pi s_0, \quad T_2:=s_0 \pi:
\end{equation}
\begin{align*}
T_1(*) \rightarrow &(-q_1-\frac{2\alpha_2+1}{2(p_1+t-2p_2^2+4q_1(q_2+q_1p_2))^2},-p_1-t+2p_2^2-4q_1(q_2+q_1p_2)),\\
&-q_2-4q_1(q_1^2+p_2),-p_2-2q_1^2,t;\alpha_2+1)
\end{align*}
and
\begin{align*}
T_2(*) \rightarrow &(-q_1+\frac{2\alpha_2-1}{2p_1},\\
&-p_1-t+2p_2^2-\frac{(2q_1p_1-2\alpha_2+1)(p_2+2q_1p_1p_2+2p_1q_2-2\alpha_2 p_2)}{p_1^2},\\
&-q_2-\frac{2(2q_1p_1+1-2\alpha_2)p_2}{p_1}-\frac{(2q_1p_1+1-2\alpha_2)^3}{2p_1^3},\\
&-p_2-\frac{(2q_1p_1+1-2\alpha_2)^2}{2p_1^2},t;\alpha_2-1).
\end{align*}
These translation operators act on $\alpha_2$ as follows$:$
\begin{equation}
T_1(\alpha_2)=\alpha_2+1, \quad T_2(\alpha_2)=\alpha_2-1.
\end{equation}
\end{proposition}

Finally, we study a solution of the system \eqref{eq:3} which is written by the use of known functions.

By the transformation $\pi$, the fixed solution is derived from
\begin{align}
\begin{split}
&\alpha_2=-\alpha_2,\\
&q_1=-q_1, \quad p_1=-(p_1+t-2p_2^2+4q_1(q_2+q_1p_2)),\\
&q_2=-(q_2+4q_1(q_1^2+p_2)), \quad p_2=-(p_2+2q_1^2).
\end{split}
\end{align}
Then we obtain
\begin{equation}
(q_1,p_1,q_2,p_2;\alpha_2)=\left(0,-\frac{t}{2},0,0;0\right)
\end{equation}
as a seed solution.

Applying the B{\"a}cklund transformations $T_1,T_2$, we can obtain an infinite series of the rational solutions:
\begin{center}
\begin{tabular}{|c||c|c|c|c|} \hline 
$\alpha_2$  & $q_1$ & $p_1$ & $q_2$ & $p_2$  \\ \hline 
-3 & $\frac{3(t^5+96)}{t(t^5-144)}$ & $-\frac{t(t^{10}-1008t^5-48384)}{2(t^5-144)^2}$ & $\frac{24(t^{15}+2088t^{10}+114048t^5-497664)}{t^3(t^5-144)^3}$ & $-\frac{12(t^{10}+432t^5+3456)}{t^2(t^5-144)^2}$   \\ \hline 
-2 & $\frac{2}{t}$ & $\frac{72}{t^4}-\frac{t}{2}$ & $\frac{12}{t^3}$ & $-\frac{6}{t^2}$   \\ \hline 
-1 & $\frac{1}{t}$ & $-\frac{t}{2}$ & $\frac{4}{t^3}$ & $-\frac{2}{t^2}$   \\ \hline 
0 & $0$ & $-\frac{t}{2}$ & $0$ & $0$   \\ \hline 
1 & $-\frac{1}{t}$ & $-\frac{t}{2}$ & $0$ & $0$   \\ \hline 
2 & $-\frac{2}{t}$ & $-\frac{t}{2}$ & $\frac{4}{t^3}$ & $-\frac{2}{t^2}$   \\ \hline 
3 & $-\frac{3(t^5+96)}{t(t^5-144)}$ & $-\frac{t^5-144}{2t^4}$ & $\frac{12}{t^3}$ & $-\frac{6}{t^2}$   \\ \hline 
\end{tabular}
\end{center}

The system \eqref{eq:3} with special parameter $\alpha_2=\frac{1}{2}$ admits a particular solution expressed in terms of the Painlev\'e I function $(q_2,p_2)$:
\begin{equation}
p_1=0,
\end{equation}
and
\begin{equation}
  \left\{
  \begin{aligned}
   \frac{dq_1}{dt} =&q_1^2+p_2,\\
   \frac{dq_2}{dt} =&-3p_2^2+\frac{t}{2},\\
   \frac{dp_2}{dt} =&q_2.
   \end{aligned}
  \right. 
\end{equation}

\section{Other polynomial Hamiltonian system}
In this section, we study some Hamiltonians transformed by birational and symplectic transformations $r_1$ and $r_2$ \rm{(see \eqref{holo1}) \rm} for the Hamiltonian \eqref{eq:4}, respectively.

At first, we study the following Hamiltonian system explicitly given by
\begin{equation}\label{eq:10}
  \left\{
  \begin{aligned}
   \frac{dq_1}{dt} =&\frac{\partial \tilde{H}}{\partial p_1}=-q_1^2p_2-1,\\
   \frac{dp_1}{dt} =&-\frac{\partial \tilde{H}}{\partial q_1}=2q_1p_1p_2+\left(\frac{1}{2}-\alpha_2 \right)p_2,\\
   \frac{dq_2}{dt} =&\frac{\partial \tilde{H}}{\partial p_2}=-3p_2^2+\frac{t}{2}-q_1^2p_1-\left(\frac{1}{2}-\alpha_2 \right)q_1,\\
   \frac{dp_2}{dt} =&-\frac{\partial \tilde{H}}{\partial q_2}=q_2
   \end{aligned}
  \right. 
\end{equation}
with the polynomial Hamiltonian
\begin{align}\label{eq:11}
\begin{split}
\tilde{H}=&-p_1+H_I(q_2,p_2,t)-\left(q_1p_1+\frac{1-2\alpha_2}{2} \right)q_1p_2\\
=&-p_1-\frac{q_2^2}{2}-p_2^3+\frac{t}{2}p_2-\left(q_1p_1+\frac{1-2\alpha_2}{2} \right)q_1p_2,
\end{split}
\end{align}
where $H_I(q,p,t):=-\frac{q^2}{2}-p^3+\frac{t}{2}p$ is the Painlev\'e I Hamiltonian.

\begin{proposition}
The birational and symplectic transformation $r_1$ \rm{(see \eqref{holo1}) \rm} takes the system \eqref{eq:3} into the Hamiltonian system \eqref{eq:10}.
\end{proposition}

\begin{theorem}
Let us consider a polynomial Hamiltonian system with Hamiltonian $H \in {\Bbb C}(t)[q_1,p_1,q_2,p_2]$. We assume that

$(B1)$ $deg(H)=5$ with respect to $q_1,p_1,q_2,p_2$.

$(B2)$ This system becomes again a polynomial Hamiltonian system in each coordinate system $\tilde{r_i} \ (i=0,1)${\rm : \rm}
\begin{align*}
\begin{split}
\tilde{r}_1:(x_1,y_1,z_1,w_1)=&\left(\frac{1}{q_1},-\left(q_1p_1-\frac{1}{2}+\alpha_2 \right)q_1,q_2,p_2 \right),\\
\tilde{r}_2:(x_2,y_2,z_2,w_2)=&\left(q_1,p_1+\frac{2\alpha_2}{q_1}+\frac{2p_2^2-t}{q_1^2}-\frac{4q_2}{q_1^3}-\frac{4p_2}{q_1^4},q_2+\frac{4p_2}{q_1}+\frac{4}{q_1^3},p_2+\frac{2}{q_1^2} \right).
\end{split}
\end{align*}
Then such a system coincides with the system \eqref{eq:10} with the polynomial Hamiltonian \eqref{eq:11}.
\end{theorem}
We note that the condition $(B2)$ should be read that
\begin{align*}
&\tilde{r}_1(H), \quad \tilde{r}_2 \left(H-\frac{1}{q_1} \right)
\end{align*}
are polynomials with respect to $x_i,y_i,z_i,w_i$.

\begin{theorem}
The system \eqref{eq:10} is invariant under the following transformations$:$ with {\it the notation} $(*):=(q_1,p_1,q_2,p_2,t;\alpha_2)$\rm{: \rm}
\begin{align*}
\begin{split}
s_0:(*) \rightarrow &\left(q_1+\frac{\alpha_2-\frac{1}{2}}{p_1},p_1,q_2,p_2,t;1-\alpha_2 \right),\\
s_1:(*) \rightarrow &\left(-q_1,-p_1-\frac{2\alpha_2}{q_1}-\frac{2p_2^2-t}{q_1^2}+\frac{4q_2}{q_1^3}+\frac{4p_2}{q_1^4},-q_2-\frac{4p_2}{q_1}-\frac{4}{q_1^3},-p_2-\frac{2}{q_1^2},t;-\alpha_2 \right).
\end{split}
\end{align*}
\end{theorem}

We remark that the birational and symplectic transformation $r_2$ \rm{(see \eqref{holo1}) \rm} takes the system \eqref{eq:3} into a polynomial Hamiltonian system. We can also obtain the same results on this Hamiltonian.

\section{Fifth-order ordinary differential equation satisfied the Hamiltonian}
In this section, we study a 1-parameter family of the fifth-order ordinary differential equation:
\begin{align}\label{eq:5}
\begin{split}
\frac{d^5u}{dt^5}=&((1-\alpha_2)\alpha_2-2\frac{du}{dt}\left(24\left(\frac{du}{dt}\right)^2-t \right)^2-24\frac{d^2u}{dt^2}\left(\frac{du}{dt}-t\frac{d^2u}{dt^2} \right)\\
&+8\frac{d^3u}{dt^3}\left(5\frac{du}{dt}\left(t-24\left(\frac{du}{dt}\right)^2 \right)-12\left(\frac{d^2u}{dt^2} \right)^2-16\frac{du}{dt}\frac{d^3u}{dt^3}\right)\\
&+2\frac{d^4u}{dt^4}\left(48\frac{du}{dt}\frac{d^2u}{dt^2}+2\frac{d^4u}{dt^4}-1 \right) )/\left(48\left(\frac{du}{dt} \right)^2+8\frac{d^3u}{dt^3}-2t \right).
\end{split}
\end{align}
The Hamiltonian \eqref{eq:4} $u:=H$ satisfies the equation \eqref{eq:5}.

In order to obtain the symmetry and holomorphy conditions for the equation \eqref{eq:5}, by making birational transformations we transform the system of rational type into the system of the first-order ordinary differential equations of polynomial type in dimension five:
\begin{equation}\label{eq:6}
  \left\{
  \begin{aligned}
   \frac{dx}{dt} =&y,\\
   \frac{dy}{dt} =&z,\\
   \frac{dz}{dt} =&-6y^2+w+\frac{t}{4},\\
   \frac{dw}{dt} =&wq+\frac{2\alpha_2-1}{4},\\
   \frac{dq}{dt} =&-\frac{1}{2}q^2-4y.
   \end{aligned}
  \right. 
\end{equation}
This system is a Riccati extension of the system \eqref{eq:3} by a scale transformation.

\section{Differential system of Polynomial type}

At first, we make the equation \eqref{eq:5}.
After differentiating once, we obtain
\begin{equation*}
\frac{dH}{dt}=\frac{w}{2}.
\end{equation*}
We can express $w$ by using the variables $\frac{dH}{dt}$:
\begin{equation*}
w=2\frac{dH}{dt}.
\end{equation*}
After differentiating again, we obtain
\begin{equation*}
\frac{d^2H}{dt^2}=\frac{z}{2}.
\end{equation*}
We can express $z$ by using the variables $\frac{d^2H}{dt^2}$:
\begin{equation*}
z=2\frac{d^2H}{dt^2}.
\end{equation*}
After differentiating again, we obtain
\begin{equation*}
\frac{d^3H}{dt^3}=\frac{1}{4}\left(t+2y-24\left(\frac{dH}{dt}\right)^2 \right).
\end{equation*}
We can express $y$ by using the variables $\frac{dH}{dt},\frac{d^3H}{dt^3}$:
\begin{equation*}
y=\frac{1}{2}\left(-t+24\left(\frac{dH}{dt}\right)^2+4\frac{d^3H}{dt^3} \right).
\end{equation*}
After differentiating again, we obtain
\begin{equation*}
\frac{d^4H}{dt^4}=\frac{1}{2}\left(tx+\alpha_2-24x\left(\frac{dH}{dt}\right)^2-24\frac{dH}{dt}\frac{d^2H}{dt^2}-4x\frac{d^3H}{dt^3} \right).
\end{equation*}
We can express $x$ by using the variables $\frac{dH}{dt},\frac{d^2H}{dt^2},\frac{d^3H}{dt^3},\frac{d^4H}{dt^4}$:
\begin{equation*}
x=\frac{2\frac{d^4H}{dt^4}+24\frac{dH}{dt}\frac{d^2H}{dt^2}-\alpha_2}{t-24\left(\frac{dH}{dt}\right)^2-4\frac{d^3H}{dt^3}}.
\end{equation*}
After differentiating, we obtain the equation \eqref{eq:5}.

Now, let us transform the equation \eqref{eq:5} into a system of polynomial type by birational transformations.
\begin{theorem}\label{th:9.1}
The birational transformations
\begin{equation*}
  \left\{
  \begin{aligned}
   x =&u,\\
   y =&\frac{du}{dt},\\
   z =&\frac{d^2u}{dt^2},\\
   w =&\frac{d^3u}{dt^3}-\frac{1}{4}\left(t-24\left(\frac{du}{dt}\right)^2 \right),\\
   q =&\frac{\frac{d^4u}{dt^4}-\frac{\alpha_2-24\frac{du}{dt}\frac{d^2u}{dt^2}}{2}}{\frac{d^3u}{dt^3}-\frac{1}{4}\left(t-24\left(\frac{dH_{III}}{dt}\right)^2 \right)}
   \end{aligned}
  \right. 
\end{equation*}
take the system of rational type into the system \eqref{eq:6}.
\end{theorem}

{\bf Proof.} At first, we rewrite the equation \eqref{eq:5} to the system of the first-order ordinary differential equations.

{\bf Step 0:} We make a change of variables.
\begin{equation}
x=u, \quad y=\frac{du}{dt}, \quad z=\frac{d^2u}{dt^2}, \quad w=\frac{d^3u}{dt^3}, \quad q=\frac{d^4u}{dt^4}.
\end{equation}

{\bf Step 1:} We make a change of variables.
\begin{equation}
x_1=x, \quad y_1=y, \quad z_1=z, \quad w_1=w+6y^2-\frac{t}{4}, \quad q_1=q.
\end{equation}
In this coordinate system, we see that this system has two accessible singular loci:
\begin{equation}
(x_1,y_1,z_1,w_1,q_1)=\left\{\left(x_1,y_1,z_1,w_1,-12y_1z_1+\frac{1}{2}\alpha_2 \right),\left(x_1,y_1,z_1,w_1,-12y_1z_1+\frac{1-\alpha_2}{2} \right) \right\}.
\end{equation}

{\bf Step 2:} We make a change of variables.
\begin{equation}
x_2=x_1, \quad y_2=y_1, \quad z_2=z_1, \quad w_2=w_1, \quad q_2=q_1+12y_1z_1-\frac{1}{2}\alpha_2.
\end{equation}
In this coordinate system, we can rewrite the system satisfying the condition \eqref{b}:
\begin{align*}
\frac{d}{dt}\begin{pmatrix}
             x_2 \\
             y_2 \\
             z_2 \\
             w_2 \\
             q_2
             \end{pmatrix}&=\frac{1}{w_2}\left\{\begin{pmatrix}
             0 & 0 & 0 & 0 & 0  \\
             0 & 0 & 0 & 0 & 0 \\
             0 & 0 & 0 & \frac{t}{4} & 0 \\
             0 & 0 & 0 & \frac{\alpha_2}{2}-\frac{1}{4} & 0 \\
             0 & 0 & 0 & 0 & \frac{\alpha_2}{2}-\frac{1}{4}
             \end{pmatrix}\begin{pmatrix}
             x_2 \\
             y_2 \\
             z_2 \\
             w_2 \\
             q_2
             \end{pmatrix}+\cdots\right\},
             \end{align*}
and we can obtain the local index $\left(0,0,0,\frac{\alpha_2}{2}-\frac{1}{4},\frac{\alpha_2}{2}-\frac{1}{4} \right)$ at the point $\{(x_2,y_2,z_2,w_2,q_2)=(0,0,0,0,0)\}$. The continued ratio of the local index at the point $\{(x_2,y_2,z_2,w_2,q_2)=(0,0,0,0,0)\}$ are all positive integers
\begin{equation}
\left(\frac{0}{\frac{\alpha_2}{2}-\frac{1}{4}},\frac{0}{\frac{\alpha_2}{2}-\frac{1}{4}},\frac{0}{\frac{\alpha_2}{2}-\frac{1}{4}},\frac{\frac{\alpha_2}{2}-\frac{1}{4}}{\frac{\alpha_2}{2}-\frac{1}{4}} \right)=(0,0,0,1).
\end{equation}

Let us resolve this accessible singular locus.

{\bf Step 3:} We blow up along the 3-fold $\{(x_2,y_2,z_2,w_2,q_2)=(x_2,y_2,z_2,0,0)\}$.
\begin{equation}
x_3=x_2, \quad y_3=y_2, \quad z_3=z_2, \quad w_3=w_2, \quad q_3=\frac{q_2}{w_2},
\end{equation}
and we can obtain the system \eqref{eq:6}.

Thus, we have completed the proof of Theorem \ref{th:9.1}. \qed

We remark that we can discuss the case of
$$
(x_1,\ldots,w_1,q_1)=\left(x_1,\ldots,w_1,-12y_1z_1+\frac{1-\alpha_2}{2} \right)
$$
in the same way.

Next, we give its symmetry and holomorphy conditions.

\begin{theorem}
The system \eqref{eq:6} admits extended affine Weyl group symmetry of type $A_1^{(1)}$ as the group of its B{\"a}cklund transformations, whose generators $s_0,s_1,{\pi}$ defined as follows$:$ with {\it the notation} $(*):=(x,y,z,w,q,t;\alpha_2)$\rm{: \rm}
\begin{align*}
s_0:(*) \rightarrow &\left(x,y,z,w,q+\frac{\alpha_2-\frac{1}{2}}{w},t;1-\alpha_2 \right),\\
s_1:(*) \rightarrow &(x+\frac{2\alpha_2+1}{2(t+2w+2yq^2-8y^2-4zq)},\\
&y+\frac{(2\alpha_2+1)q}{2(t+2w+2yq^2-8y^2-4zq)}-\frac{(2\alpha_2+1)^2}{4(t+2w+2yq^2-8y^2-4zq)^2},\\
&z+\frac{(2\alpha_2+1)(q^2-8y)}{4(t+2w+2yq^2-8y^2-4zq)}-\frac{3(2\alpha_2+1)^2 q}{4(t+2w+2yq^2-8y^2-4zq)^2}\\
&+\frac{(2\alpha_2+1)^3}{4(t+2w+2yq^2-8y^2-4zq)^3},\\
&w+\frac{2(2\alpha_2+1)(yq-z)}{t+2w+2yq^2-8y^2-4zq}+\frac{(2\alpha_2+1)^2(q^2+4y)}{4(t+2w+2yq^2-8y^2-4zq)^2},\\
&q-\frac{2\alpha_2+1}{t+2w+2yq^2-8y^2-4zq},t;-1-\alpha_2),\\
\pi:(*) &\rightarrow (x+\frac{q}{2},-\left(y+\frac{q^2}{4} \right),-\left(z-\frac{1}{4}q(q^2+8y) \right),\\
&-\left(w+yq^2-4y^2-2zq+\frac{t}{2} \right),-q,t;-\alpha_2).
\end{align*}
\end{theorem}

\begin{theorem}
Let us consider a system of first order ordinary differential equations in the polynomial class\rm{:\rm}
\begin{align*}
&\frac{dx}{dt}=f_1(x,y,z,w,q), \quad \frac{dy}{dt}=f_2(x,y,z,w,q), \ldots ,\frac{dq}{dt}=f_5(x,y,z,w,q),\\
&f_i \in {\Bbb C}(t)[x,y,z,w,q] \ (i=1,\cdots,5).
\end{align*}
We assume that

$(A1)$ $deg(f_i)=2$ with respect to $x,y,z,w,q$.

$(A2)$ The right-hand side of this system becomes again a polynomial in each coordinate system $r_i:(x_i,y_i,z_i,w_i,q_i) \ (i=1,2,3):$
\begin{align*}
r_1:(x_1,y_1,z_1,w_1,q_1)=&\left(x,y,z,-\left(wq+\frac{1}{2}(2\alpha_2-1) \right)q,\frac{1}{q} \right),\\
r_2:(x_2,y_2,z_2,w_2,q_2)=&(x+\frac{q}{2},y+\frac{q^2}{4},z-\frac{1}{4}q(q^2+8y),\\
&-\left(\left(w+\frac{t}{2}-4y^2+q(yq-2z) \right)q-\frac{1}{2}-\alpha_2 \right)q,\frac{1}{q}),
\end{align*}
\begin{align*}
r_3:(x_3,y_3,z_3,w_3,q_3)=&(\frac{1}{x},\frac{y}{x}+x,-x^3z-\frac{1}{64}y(4t+5w-35y^2)\\
&-\frac{1}{16}x(4wq+2\alpha_2+1)+\frac{3y(y^2-w)^2}{128x^4}-\frac{105}{128}x^4y-\frac{7}{128}x^6\\
&-\frac{3}{64}x^2(4t+15w-35y^2)+\frac{(w-21y^2)(w-y^2)}{128x^2},\\
&\frac{w}{x^2},qx^2+\frac{3}{4}x(x^2-2y)+\frac{w-y^2}{4x} ).
\end{align*}
Then such a system coincides with the system
\begin{equation*}
  \left\{
  \begin{aligned}
   \frac{dx}{dt} =&y+g(t) \quad (g(t) \in {\Bbb C}(t)),\\
   \frac{dy}{dt} =&z,\\
   \frac{dz}{dt} =&-6y^2+w+\frac{t}{4},\\
   \frac{dw}{dt} =&wq+\frac{2\alpha_2-1}{4},\\
   \frac{dq}{dt} =&-\frac{1}{2}q^2-4y.
   \end{aligned}
  \right. 
\end{equation*}
\end{theorem}

Finally, we study a solution of the system \eqref{eq:6} which is written by the use of known functions.

By the transformation $\pi$, the fixed solution is derived from
\begin{align}
\begin{split}
&\alpha_2=-\alpha_2,\\
&x+\frac{q}{2}=x, \quad -\left(y+\frac{q^2}{4} \right)=y, \quad -\left(z-\frac{1}{4}q(q^2+8y) \right)=z,\\
&-\left(w+yq^2-4y^2-2zq+\frac{t}{2} \right)=w, \quad -q=q.
\end{split}
\end{align}
Then we obtain
\begin{equation}
(x,y,z,w,q;\alpha_2)=\left(x,0,0,-\frac{t}{4},0;0 \right)
\end{equation}
as a seed solution.

The system \eqref{eq:6} with special parameter $\alpha_2=\frac{1}{2}$ admits a particular solution expressed in terms of the Painlev\'e I function $(y,z)$:
\begin{equation}
w=0,
\end{equation}
and
\begin{equation}
\frac{dx}{dt} =y, \quad \frac{dy}{dt} =z, \quad \frac{dz}{dt} =-6y^2+\frac{t}{4}, \quad \frac{dq}{dt} =-\frac{1}{2}q^2-4y.
\end{equation}

\end{document}